\documentclass[letterpaper, 10 pt, conference]{ieeeconf}  

\IEEEoverridecommandlockouts                              

\usepackage{graphicx} 
\usepackage{epsfig} 
\usepackage{times} 
\usepackage{amsmath} 
\usepackage{amssymb}  
\usepackage{float}
\usepackage{caption}
\usepackage{subcaption}
\usepackage{algorithm}
\usepackage{algpseudocode}
\usepackage{xcolor} 

\title{\LARGE \bf Optimal Ensemble Control of Matter-Wave Splitting \\ in Bose-Einstein Condensates
}

\author{Andre Luiz P. de Lima$^1$, Andrew K. Harter$^2$, Michael J. Martin$^2$, and Anatoly Zlotnik$^3$
\thanks{This project was supported by the LDRD program and the Center for Nonlinear Studies at Los Alamos National Laboratory.  Research conducted at Los Alamos National Laboratory is done under the auspices of the National Nuclear Security Administration of the U.S. Department of Energy under Contract No. 89233218CNA000001.}
\thanks{$^{1}$Andre Luiz P. de Lima is a Ph.D. student in the Electrical \& Systems Engineering Department at Washington University in St. Louis, St. Louis, MO, USA
    {\tt\small a.delima@wustl.edu}}%
\thanks{$^{2}$Andrew K. Harter and Michael J. Martin are with the Quantum Technologies Team at Los Alamos National Laboratory, Los Alamos, NM, USA
    {\tt\small \{harter, mmartin\}@lanl.gov}}%
\thanks{$^3$Anatoly Zlotnik is in the Applied Mathematics \& Plasma Physics Group, Los Alamos National Laboratory, Los Alamos, NM, USA
        {\tt\small azlotnik@lanl.gov}}%
        }

\begin{document}

\maketitle
\thispagestyle{empty}
\pagestyle{empty}

\begin{abstract}
    We present a framework for designing optimal optical pulses for the matter-wave splitting of a Bose-Einstein Condensate (BEC) under the influence of experimental inhomogeneities, so that the sample is transferred from an initial rest position into a singular higher diffraction order. To represent the evolution of the population of atoms, the Schrödinger's equation is reinterpreted as a parameterized ensemble of dynamical units that are disparately impacted by the beam light-shift potential in a continuous manner. The derived infinite-dimensional coupled Raman-Nath equations are truncated to a finite system of diffraction levels, and we suppose that the parameter that defines the inhomogeneity in the control applied to the ensemble system is restricted to a compact interval. We first design baseline square pulse sequences for the excitation of BEC beam-splitter states following a previous study, subject to dynamic constraints for either a nominal system assuming no inhomogeneity or for several samples of the uncertain parameter. We then approximate the continuum state-space of the ensemble of dynamics using a spectral approach based on Legendre moments, which is truncated at a finite order.  Control functions that steer the BEC system from an equivalent rest position to a desired final excitation are designed using a constrained optimal control approach developed for handling nonlinear dynamics. This representation results in a minimal dimension of the computational problem and is shown to be highly robust to inhomogeneity in comparison to the baseline approach.  Our method accomplishes the BEC-splitting state transfer for each subsystem in the ensemble, and is promising for precise excitation in experimental settings where robustness to environmental and intrinsic noise is paramount.
\end{abstract}


\section{INTRODUCTION} \label{sec:intro}



Numerous challenges in quantum science and technology involve design, observation, or control of bilinear dynamical systems \cite{koch2022quantum}. A compelling application in metrology involves cold atom interferometry, in which optical standing waves are used to split and recombine matter waves \cite{Cassidy:21}.  In this setting, a time-varying laser pulse is applied to steer a quantum system between states of interest.  Such experiments demonstrate significant sensitivity to accelerations and rotations, and have the potential to provide real improvement in measurements of fundamental physics \cite{canuel2006six}.  However, experimental uncertainties are inherent in these systems, and this reduces the precision in the splitting of matter waves, which can lead to inconsistent diffraction fidelity \cite{Altin:13}.


Robustness to experimental inhomogeneities, such as laser intensity, was investigated in atom interferometry. Previous approaches have attempted to minimize the influence of these dynamical disturbances, using for example adiabatic transfer techniques, which do not depend on exact optical intensity \cite{Jaffe:18}. Other methods seek to optimally design light pulses to improve tolerance to systematic inhomogeneities by using specific control profile shapes, such as composite pulses \cite{Dunning:14}, or large-scale optimization, such as GRadient Ascent Pulse Engineering (GRAPE) \cite{Saywell:20}. A practical means to parametrize  inhomogeneities is by way of a factor that scales the light amplitude control applied to the system dynamics. This mathematical description enables the state evolution of the sample to be represented by a collection of ordinary differential equations, each of which describes the possible evolution taken by an individual atom in the experiment subject to the inexact effect of the applied optical pulse. 


Ensemble control theory was developed to address the need to understand the evolution of systems comprised by a very large (potentially infinite) number of similar dynamical units, such as molecular, atomic and quantum systems, when manipulated by a single universally applied control field \cite{Li:06}. This area of mathematical systems theory focuses on defining controllability conditions \cite{Zeng:16}, measurement processing \cite{Yu:20}, and control design approaches \cite{Leghtas:11} for parameterized dynamical systems. Ensemble control theory has been used to address challenging control tasks for large-scale dynamical problems, including brain dynamics \cite{Zlotnik:16} and Magnetic Resonance Imaging (MRI) \cite{Li:09}. Furthermore, ensemble control was shown to be effective in pulse design for Nuclear Magnetic Resonance (NMR) accounting for the influence of systematic disturbances \cite{Andre:22, Li:09}.  The mathematical setting used for control synthesis for bilinear ensemble systems \cite{Li:09} can be applied to develop laser pulses for the atom interferometry procedure, and indeed the simplification of the Raman-Nath equations (RNE) truncated to a two-level system results in the Bloch equations in NMR \cite{Wu:05,zlotnik2012iterative}.


Ensemble systems are challenging because of infinite-dimensionality arising from a parameter that varies continuously on a compact interval. Whereas the state of each individual ensemble element is finite-dimensional, the state of an ensemble system lies on a continuous Hilbert space. In this context, methods designed for unified representations were studied with the goal of reducing computational complexity, and polynomial moments can be used to significantly reduce the size of representations of ensemble system dynamics \cite{Narayanan:20}.  Moment dynamics entail spectral approximation of the ensemble in Hilbert space, with the spectral accuracy property that ensures exponential convergence to the continuum state with increasing approximation order.  This representation enables ensemble control design in a space of greatly reduced dimension without significant loss of precision \cite{Li:22}. 


In this paper, we revisit the optimal control problem (OCP) of matter-wave splitting involving a Bose-Einstein condensate (BEC) in a standing light wave potential \cite{Wu:05}.  We consider optical standing-wave beam splitters for which the effect of an optical pulse is subject to inhomogeneities, which we represent by an uncertain factor applied to the light wave amplitude envelope.  We describe the evolution of the resulting dynamic ensemble using Legendre moment dynamics, then design pulse sequences that transfer the state of the ensemble to a desired set-point, and examine the improved fidelity of the BEC momentum transfer within the quantum state space.  We demonstrate that our new pulse design approach can withstand up to $\pm10$\% perturbation to the applied control field applied to a BEC interferometry experiment for the examined diffraction orders, and offers unprecedented flexibility with respect to previous studies \cite{Cassidy:21}.


This project is organized as follows.  In Section \ref{sec:system}, we derive the mathematical representation of the BEC-splitting system from the Schrödinger's equation, and characterize the source of inhomogeneity. In Section \ref{sec:square}, we then reproduce a state-of-the-art method of using square pulse sequences for optical splitting of matter-waves \cite{Cassidy:21}, and extend it to compensate for parameter inhomogeneity using sampling. In Section \ref{sec:moments}, we describe the Legendre moment representation of the BEC splitting dynamics, our proposed method for light pulse design within the moment space, and the resulting optimal control formulation. We present the computed optimal controls and simulations of the resulting momentum transfers in Section \ref{sec:results} for several diffraction orders, which show a significant improvement in fidelity with respect to the benchmark method.  Finally, we conclude in Section \ref{sec:conc}.


\section{ENSEMBLE SYSTEM FOR BEC SPLITTING} \label{sec:system}


We consider a dilute BEC composed of an initially stationary population of atoms, which is exposed to a standing-wave light field that excites the cluster of particles with similar, yet distinct intensities. The governing dynamics for this system are described by a one-dimensional Schrödinger's equation for a single atom. Inhomogeneity in the light beam is represented by a multiplicative parameter $\varepsilon$ that scales the light shift potential amplitude $\Omega(t)$, which leads to a wave function
\begin{equation}\label{eq:schrodinger}
   \!\! i\dot{\psi}(x, t, \varepsilon) \!=\! \left( -\frac{\hbar}{2m}\frac{d^{2}}{dx^{2}} \!+\! \varepsilon\Omega(t)cos\left(2k_{0}x\right) \right)\psi(x, t, \varepsilon).
\end{equation}
The time-dependent parameter $\Omega(t)$ serves as the control input available to the experimenter, and $k_{0}$ is the wave vector of the light field. We suppose that $\varepsilon \in K$ is the physical representation of the unequal influence of the control pulse on individual atoms in the BEC sample, where $K \equiv [1-\delta,1+\delta]$ is a compact interval. Each atomic particle can be viewed as corresponding to a unique $\varepsilon$, so that this parameter can be used as an index for the dynamical units in the sample. 

Equation \eqref{eq:schrodinger} can be further developed into a sequence of coupled Raman-Nath equations by describing the wave function within the Bloch basis. The resulting differential equations dictate the time evolution of the beam-splitter state populations (represented by ground state $C_{0}$, and symmetric and anti-symmetric superpositions of momentum states $C_{2n}^{+}$ and $C_{2n}^{-}$,  for $n=0,1,2,\ldots$) and the manner in which these interact.  This infinite-dimensional system of ordinary differential equations (ODEs) is \cite{Wu:05}
\begin{multline}\label{eq:Raman-Nath}
    i\dot{C_{2n}}(k, t) = \frac{\hbar}{2m}(2nk_{0}+k)^{2}C_{2n}(k, t)\\ + \frac{\varepsilon\Omega(t)}{2}[C_{2n-2}(k, t)+C_{2n+2}(k, t)].
\end{multline}
We make several assumptions in order to further simplify the infinite-dimensional Raman-Nath equation system:
\begin{itemize}
    \item The dependence on parameter $k$ is dropped by assuming the initial momentum to be narrowly distributed around the value of $k=0$.
    \item The coupling of the states $C_{0}$ and $C_{2n}^{+}$ with their anti-symmetric counterparts $C_{2n}^{-}$ is considered negligible by assuming that $k \ll k_{0}$.
    \item A value $N^{+}$ is defined such that all levels above $C_{2N^{+}}^{+}$ are permanently unpopulated. This is achievable by assuming that all atoms are at rest at the beginning and that $\Omega(t)/2 \ll (2N^{+})^{2}\hbar k_{0}^{2}/2m = (2N^{+})^{2} \omega_{r}$, where $\omega_{r}$ is the photon recoil frequency.
\end{itemize}
This mathematical setting has been developed in a variety of previous studies on matter-wave splitting techniques \cite{cronin2009optics,Wu:05,Cassidy:21}. The resulting state-space system is described by

\begin{equation}\label{eq:BEC_ensemble}
    \frac{d}{dt}\begin{bmatrix}
    C_{0}(t, \varepsilon)\\C_{2}^{+}(t, \varepsilon)\\ \vdots \\ C_{2N}^{+}(t, \varepsilon)
    \end{bmatrix} = -i\omega_{r}A(\varepsilon, \Omega(t)) 
    \begin{bmatrix}
    C_{0}(t, \varepsilon)\\C_{2}^{+}(t, \varepsilon)\\ \vdots \\ C_{2N}^{+}(t, \varepsilon)
    \end{bmatrix},
\end{equation}
where $A(t, \varepsilon)$ is a real $(N^{+}+1 \times N^{+}+1)$ symmetric matrix defined by
\begin{multline} \label{eq:BEC_ensemble_A}
    A(\varepsilon, \Omega(t)) =     \begin{bmatrix}
    0& 0& 0& \hdots& 0\\0 & 4 & 0 & \hdots & 0\\ 0 & 0 & 16 & \hdots & 0 \\ \vdots & \vdots & \vdots & \ddots & \vdots \\ 0 & 0 & 0 & \hdots & (2N^{+})^{2}
    \end{bmatrix}\\ + \frac{\varepsilon\Omega(t)}{2\omega_{r}}
    \begin{bmatrix}
    0& \sqrt{2}& 0& 0& \hdots& 0\\ \sqrt{2} & 0 & 1 & 0& \hdots & 0\\ 0 & 1 & 0 & 1 & \hdots & 0 \\0 & 0 & 1 & 0 & \ddots & \vdots \\ \vdots & \vdots & \vdots & \ddots & \ddots & 1 \\ 0 & 0 & 0 & \hdots & 1 & 0
    \end{bmatrix}.
\end{multline}

For the system in Equations \eqref{eq:BEC_ensemble}-\eqref{eq:BEC_ensemble_A}, we define the state vector $C^{+} = [C_{0}, C_{2}^{+}, ..., C_{2N^{+}}^{+}] \in \mathbb{C}^{N^{+}+1}$.  A goal of this study is to enable constrained optimization approaches for pulse design, for which the complex dynamics in Equations \eqref{eq:BEC_ensemble}-\eqref{eq:BEC_ensemble_A} must be converted to a real-valued system. We thus define $C^{+}_{\mathbb{R}} \in \mathbb{R}^{2(N^{+}+1)}$ such that $C^{+}_{\mathbb{R}} = [C_{0, Re}, C_{0, Im}, C^{+}_{2, Re}, C^{+}_{2, Im}, ..., C^{+}_{2N^{+}, Re}, C^{+}_{2N^{+}, Im}]$, in terms of which the time evolution dynamics are

\begin{equation}\label{eq:BEC_ensemble_Real}
    \frac{d}{dt}C^{+}_{\mathbb{R}} = \omega_{r} A(\varepsilon, \Omega(t)) \otimes \begin{bmatrix}
        0 & -1\\ 1 & 0
    \end{bmatrix} C^{+}_{\mathbb{R}},
\end{equation}
where $\otimes$ is the Kronecker product.

In the above representation, the momentum dynamics for each atom are truncated to a finite-dimensional quantum system of $N^{+}+1$ levels. Nonetheless, due to the presence of parameter $\varepsilon$, there is a collection of infinitely many systems that must be steered simultaneously by a common control input.  Next, we examine methods to design control pulses for momentum transfer of the system \eqref{eq:BEC_ensemble_Real} that are invariant to a continuum of values of the uncertain parameter $\varepsilon$.



\section{SPLITTING MATTER USING SQUARE PULSES}  \label{sec:square}


A previously proposed approach for splitting matter waves uses a sequence of two square pulses for the envelope of an optical input that is applied to a dilute BEC in a standing light wave potential.  The approach has been shown to have a high degree of precision in simulation, and can be applied in practice \cite{Wu:05, Xiong:11}.  The drawback of this method is its use of a single nominal system to define the dynamic response of each atom in the BEC to the optical pulse, so there is no compensation for inherent inhomogeneities. We describe the square-pulse design procedure below, and extend the original methodology with pulses of equal amplitude to the design of sequences of pulses with different amplitudes, which has been shown to enable higher fidelity \cite{Cassidy:21}.

\subsection{Square Optical Pulse Design for a Single System}   \label{sec:singlesq}

\begin{figure}[!h]
    \centering
    \includegraphics[width=0.8\linewidth]{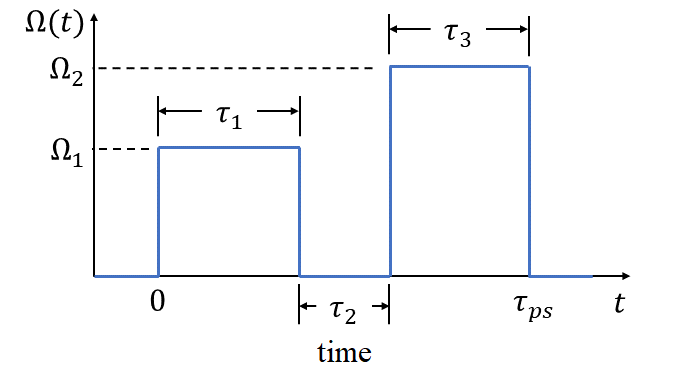}
    \vspace{-1ex}
    \caption{Illustration of a square pulse sequence control $\Omega(t)$ and parameters.}
    \label{fig:Square_Pulse_ex}
\end{figure}

The square pulse optimal design problem involves the optimization of five parameters, which are illustrated in Figure \ref{fig:Square_Pulse_ex}. The decision variables are the time durations $\tau_{1}$ and $\tau_{3}$ of the pulses,  the time interval $\tau_{2}$ between the pulses, and the amplitudes $\Omega_{1}$ and $\Omega_{2}$ of the first and second square pulses, respectively. These parameters determine the final state of the system after the pulse is employed, which can be described by continuously applying the variation of parameters formula to Equation \eqref{eq:BEC_ensemble}.  The state at the terminal time $\tau_{ps}=\tau_1+\tau_2+\tau_3$ of the pulse sequence is
\begin{equation}\label{eq:time_evolution}
\begin{gathered}
    C^{+}(\tau_{ps}, \varepsilon) = e^{-i\omega_{r}\tau_{3}A_{3}}e^{-i\omega_{r}\tau_{2}A_{2}}e^{-i\omega_{r}\tau_{1}A_{1}}C^{+}(0, \varepsilon)\\
    \text{s.t. } A_{1} = A(\varepsilon, \Omega_{1}), \, A_{2} =A(\varepsilon, 0),\\
    A_{3} =A(\varepsilon, \Omega_{2}), \, \tau_{ps} = \tau_{1} + \tau_{2} + \tau_{3}.
\end{gathered}
\end{equation}

In the baseline scenario, we suppose that $\varepsilon = 1$.  The control design goal is to transfer the momentum state of the sample from an initial rest state $C^{+}(0) = [1, 0, 0, \hdots, 0]^{T}$ to a final desired state $C_{f}^{+}$. This desired state is typically defined with respect to the $L^{2}$-norm (in quantum momentum space), in which an energy level $2n$ is reached at the end of the sequence of pulses, i.e., $C_{f}^+=[C_{f,0}^+,C_{f,2}^+,\ldots,C_{2N^+}^+]$ with $\lVert C_{f}^{+} \rVert_{2} = \lVert C_{f, 2n}^{+} \rVert_{2} = 1$ and $\lVert C_{f,2k}^{+} \rVert_{2} = 0$ for $k\neq n$. For this control task, we formulate the optimization problem as

\begin{equation} \label{eq:Square_Design_one}
\begin{aligned}
    \min_{\Omega_{1}, \Omega_{2}, \tau_{1}, \tau_{2}, \tau_{3}} \quad & \lVert |C_{f}^{+}| - |C^{+}(\tau_{ps}, 1)| \rVert_{2} \text{ as in Eq.} \eqref{eq:time_evolution} \\
    s.t. \quad & \textrm{Dynamics in Equations \eqref{eq:BEC_ensemble_A}-\eqref{eq:BEC_ensemble_Real}}.
\end{aligned}
\end{equation}


The procedure described above follows the original methodology for the design of square pulse sequences for matter wave splitting. Because a nominal parameter $\varepsilon$ is used, the approach does not account for non-uniform influence of the pulse on the dynamical behavior of individual atoms in the BEC, and we expect that control performance is not robust to values of $\varepsilon$ away from unity. We extend the formulation below by adding dynamic constraints for samples of $\varepsilon\subset K$ as an attempt to account for the dynamics of an ensemble of atoms that are affected inhomogeneously by the applied pulse sequence.

\vspace{-2ex}
\subsection{Square Optical Pulse Design using Ensemble Sampling}  \label{sec:multisq}

A widely used approach to compensating for parameter dispersion in quantum ensembles involves sampling the ensemble space \cite{chen2014sampling}.  We extend the optimization problem \eqref{eq:Square_Design_one} defined in the previous section to account for multiple values of $\varepsilon$. Because it is intractable to define dynamic constraints for each system in an infinite-dimensional ensemble, we define an objective function as a summation of minimum $L_2$ norm objectives for $m$ values of $\varepsilon$ sampled directly from the parameter domain $K$. The optimization is conducted subject to dynamic constraints for each of the same $m$ values of $\varepsilon$.

We define the finite collection of parameters $\varepsilon_{1:m} = \{\varepsilon_{1}, \varepsilon_{2}, \hdots, \varepsilon_{m}\} \subset K$, such that $\varepsilon_{1} =\mathop{\arg \min}_{K} \varepsilon$, $\varepsilon_{m} =\mathop{\arg \max}_{K} \varepsilon$, with intermediate samples uniformly spaced on $[\varepsilon_{1}, \varepsilon_{m}]$. The simultaneous optimization for all systems indexed by $i\in\{1,2,\ldots,m\}$ is formulated as 

\vspace{-1.5ex}
\begin{equation} \label{eq:Square_Design_multi}
\begin{aligned}
    \min_{\Omega_{1}, \Omega_{2}, \tau_{1}, \tau_{2}, \tau_{3}} \quad & \sum_{i=1}^{m}\lVert |C_{f}^{+}|- |C^{+}(\tau_{ps}, \varepsilon_{i})| \rVert_{2}\\
    s.t. \quad & \textrm{Dynamics in Eqns. \eqref{eq:BEC_ensemble_A}-\eqref{eq:BEC_ensemble_Real} \text{ for } $\varepsilon=\varepsilon_i$}, \\ & \qquad \forall i=1,\ldots,m.
\end{aligned}
\end{equation}
The control pulse sequences obtained by solving the optimization problems \eqref{eq:Square_Design_one} and \eqref{eq:Square_Design_multi} will be used as benchmarks in Section \ref{sec:results} to demonstrate the performance improvement gained by way of the proposed moment dynamics method. 


\section{CONSTRAINED ENSEMBLE CONTROL USING MOMENT DYNAMICS} \label{sec:moments}


We present an improved light pulse design method for matter-wave splitting by compensating for the inherent inhomogeneities present in physical systems. In this section, we derive a polynomial moment based representation of the ensemble dynamics of the matter wave splitting process, and formulate an ensemble control problem for optimal state transfer of a BEC sample in momentum space. 

The polynomial moment ensemble approach has been applied recently to design open-loop controls to transfer ensemble systems between initial and target states in Hilbert space, demonstrating low terminal state error \cite{Narayanan:20, Yu:23, Zeng:16}. The method enables a control problem for an infinite-dimensional continuum to be reduced to a finite-dimensional control problem by transforming the state space using a polynomial approximation that can then be truncated at a sufficiently high order.  The transformation facilitates the design of a single control input that steers the moments to a desired representation in the spectral moment space, which is equivalent to steering each individual element of the ensemble to its target in the state space \cite{Andre:22}.

The moment-based approach for ensemble control systems can be defined using any basis from a diverse collection of orthogonal polynomials. For this project, we have chosen to use Legendre polynomials as the basis of the moment system, because this particular basis has desirable properties for approximation of continuous functions on compact intervals \cite{Li:22}.  This representation has also been used previously for controlling the Bloch equations that represent a two-level quantum system in the context of NMR \cite{Andre:22}.

\subsection{Moment Systems using Legendre Polynomials}

The concept of moment space representation begins with a function defined within a separable Hilbert space $\mathcal{H}$. If a basis \{$|\psi_{k}\rangle: k \in \mathbb{N}^{d}$\} can be defined (for simplicity, we will use $d = 1$), then the moment of order $k$ for a function $x(t)$ can be defined by
\begin{equation} \label{eq:Moment_inner_product}
    m_{k}(t) = \langle x(t) | \psi_{k} \rangle.
\end{equation}
For an ensemble of systems indexed by a parameter varying in the compact interval $K = [-1, 1]$, it is possible to obtain the equivalent Legendre moments by using the respective orthogonal polynomials $P_{k}(\varepsilon)$, obtained using the recursive relation given, after normalization, by
\begin{equation} \label{eq:Leg_recursive}
    \varepsilon P_{k} = c_{k-1}P_{k-1} + c_{k}P_{k+1},
\end{equation}
where $P_{0}(\varepsilon) = 1/\sqrt{2}$, $P_{1}(\varepsilon) = \sqrt{3/2}\varepsilon$, $\hdots$, are the normalized Legendre polynomials and $c_{k} = (k+1)/\sqrt{(2k+1)(2k+3)}$. In this setting, the Legendre moment of order $k$ for the vector of quantum states $C^{+}(t, \varepsilon)$ is defined by
\begin{equation} \label{eq:Leg_Moment}
    m_{k}(t) = \langle C^{+}(t), P_{k} \rangle = \int_{-1}^{1} C^{+}(t, \varepsilon)P_{k}(\varepsilon) d\varepsilon.
\end{equation}
Taking the time derivative of Equation \eqref{eq:Leg_Moment} yields the equivalent moment dynamics of the ensemble system. The use of the Legendre polynomial basis to define moment spaces was shown to possess a variety of advantages for representing ensemble systems. For instance, the Legendre polynomials form an orthogonal basis on the Hilbert space $\mathcal{H}$, which induces an isometry between the moment space and the quantum ensemble system.  The preservation of metrics through this transformation has advantages for defining tractable optimization formulations to represent OCPs for ensemble systems. The recursive relation in Equation \eqref{eq:Leg_recursive} also enables the description of moment dynamics defined by bounded operators. Finally, the control profile is preserved with transformation between the moment space and state space, so that practical control design is possible using the moment dynamics method.

\subsection{BEC Splitting in the Moment Space}

We apply the transformation in Equation \eqref{eq:Leg_Moment} to the time-evolution model for the quantum system described by Equations \eqref{eq:BEC_ensemble_A}-\eqref{eq:BEC_ensemble_Real} to obtain the dynamics of the corresponding Legendre moments. The quantum ensemble system is defined for the ensemble parameter $\varepsilon \in K = [1-\delta, 1+\delta]$, such that $\delta \in (0, 1)$, which is a different compact interval than the domain of the Legendre polynomial basis. Therefore, we define a parameter $\varepsilon^{*}$ such that $\varepsilon = 1 + \delta\varepsilon^{*}$, meaning that $\varepsilon^{*} \in [-1,1]$, and the moments are calculated using the Hilbert space related to parameter $\varepsilon^{*}$. In this setting, the moment dynamics for $m(t) = [m_{0}(t), m_{1}(t), \hdots, m_{N - 1}(t)]$ are described by

\begin{align}
    &\dot{m(t)} =  \omega_{r} \left[ I(N) \otimes A(1, \Omega(t)) \otimes \begin{bmatrix}
        0 & -1\\ 1 & 0
    \end{bmatrix} \right. \label{eq:BEC_moment} \\ & \quad \quad \left. + \mathcal{C}(N) \otimes \left( A(\delta, \Omega(t)) - A(\delta, 0)\right)  \otimes \begin{bmatrix}
        0 & -1\\ 1 & 0
    \end{bmatrix} \right] \!\cdot\!  m(t), \nonumber
\end{align}
where $I(N)$ is an identity matrix of dimension $N \times N$ and $\mathcal{C}(N)$ is defined by

\begin{equation} \label{eq:BEC_moment_C}
    \mathcal{C}(N) = \begin{bmatrix}
        0 & c_{0} & & & \\
        c_{0} & 0 & c_{1} & & \\
        & c_{1} & \ddots & \ddots & &\\
        & & \ddots & 0 & c_{N-2}\\
        & & & c_{N-2} & 0
    \end{bmatrix}.
\end{equation}

With the Legendre moment dynamics as defined in Equation \eqref{eq:BEC_moment}, we can state an optimization problem that approximates the optimal control problem for the ensemble in state-space. The orthogonality property of the Legendre polynomial basis defines an isometry between the two studied spaces, which enables a metric to be used to specify an objective for the optimization problem. Moreover, it is also possible to include inequality constraints in the optimization problem, which we will apply in the control design.  We are interested specifically in enforcing limits on the absolute amplitude and the rate of change of the control function.

Specifying the final state in moment space requires an additional reformulation. In the state-space, the optimization objective is defined in terms of the absolute value of the final state, as in Equation \eqref{eq:Square_Design_one}. To express an objective function that depends on the state of an ensemble system,  we formulate the objective in moment space.  To that end, we minimize the magnitude of the undesired energy levels by aiming to nullify their states, i.e.,  $\lVert C_{f}^{+} \circ (\mathbf{1}_{2(N^{+}+1)} - e_{n+1}\otimes \mathbf{1}_{2}) \rVert_{2}  = 0$. The resulting OCP in moment space is
\begin{equation} \label{eq:moment_optimization}
\begin{aligned}
    \min_{\Omega(t)} \quad & \lVert m(T) \circ [\mathbf{1}_{N} \otimes (\mathbf{1}_{2(N^{+}+1)} - e_{n+1} \otimes \mathbf{1}_{2})] \rVert_{2}\\
    s.t. \quad & \textrm{Dynamics in Equations \eqref{eq:BEC_moment}-\eqref{eq:BEC_moment_C}},\\
    \quad & \Omega_{min} \leq \Omega(t) \leq \Omega_{max},\\
    \quad & \Delta\Omega_{min} \leq \dot{\Omega}(t)\leq \Delta\Omega_{max},
\end{aligned}
\end{equation}
where $T$ is a predefined time horizon, similar to the parameter $\tau_{ps}$ used in problems \eqref{eq:Square_Design_one} and \eqref{eq:Square_Design_multi}.


\section{CONTROL SYNTHESIS RESULTS}  \label{sec:results}


The Legendre moments method for matter-wave splitting ensemble control design is evaluated by comparing its performance with that of the square pulse sequence method. We first compare the Legendre moments approach with the controls obtained by solving the problem in Equation \eqref{eq:Square_Design_one}.  A similar comparison is then done with respect to controls obtained by solving the problem in Equation \eqref{eq:Square_Design_multi}, which seeks to account for variation in $\varepsilon$ by sampling values directly in the ensemble space $K$. To compare the performance of the various control designs in achieving the goal of transferring the ensemble to the target in Hilbert space on $K$, we define an index $I_{e}$ based on the objective function for the pulse design in the quantum state space. This performance index, which we seek to minimize, is defined as

\begin{equation} \label{eq:performance_index}
    I_{e} = \int_{1 - \delta}^{1 + \delta} \lVert |C_{f}^{+}| - |C^{+}(T, \varepsilon)| \rVert_{2} d\varepsilon.
\end{equation}

\subsection{Computational Approach}

We solve the OCP \eqref{eq:moment_optimization} using a general iterative scheme for computing optimal control inputs for nonlinear systems \cite{zeng2019iterative}.  For a general continuous-time nonlinear control system of the form $\dot{x}(t)=f(x(t),u(t))$, where $x(t)$ and $u(t)$ are state and control vectors, respectively, a zero-order hold assumption is used to define the evolution of the system using a piece-wise constant control input.  Beyond using a traditional first-order discretization and linearization, we use an approach that takes advantage of a higher-order Taylor series expansion to accurately represent linearization of nonlinear dynamics and their Jacobian at each iteration \cite{vu2020iterative}. Beginning with an initial guess of the control function, the algorithm solves an approximation of the optimization problem in \eqref{eq:moment_optimization} to yield the optimal vector $\Delta U$ of time-discretized variations in the control input that improves the objective value.  This results in a sequence of quadratic programs, which, though each can be efficiently solved, require highly complex algebraic formulation of the dynamics in Equations \eqref{eq:BEC_moment}-\eqref{eq:BEC_moment_C} when truncated at a high momentum level $N^+$ and sufficiently high polynomial moment order to represent the quantum and ensemble dynamics in their respective Hilbert spaces.  In contrast to previous applications of this approach to the control of quantum ensemble systems \cite{Li:22}, the size and complexity of the dynamics examined here requires the iterative use of symbolic algebra to compute the Taylor expansion of the nonlinear dynamics, which becomes complicated.

Once the moment dynamics and the respective symbolic expressions are obtained, the pulse is designed by the iterative application of quadratic programming to solve the optimization problem in Equation \eqref{eq:moment_optimization}. To apply the methodology defined in \cite{vu2020iterative}, the time domain $[0,T]$ is finely discretized by using a small time step $\Delta t$ into $m$ discrete values. Starting with an initial guess $U\in\mathbb{R}^m$, each iteration optimizes the variation $\Delta U \in \mathbb{R}^m$ that is then added to the control function estimate.  

\vspace{-1ex}
\begin{algorithm}[!h]
\caption{Iterative Moment Pulse Design Algorithm}\label{alg:Iterative}
\begin{algorithmic}[1]
\Require $n$, Initial pulse guess $U=U_0$, Moment dynamics $M(t, m(t), u(t))$, number of time steps $N_{steps}$, $\lambda$, $f(m(t)) = m(t) \circ [\mathbf{1}_{N} \otimes (\mathbf{1}_{2(N^{+}+1)} - e_{n+1} \otimes \mathbf{1}_{2})] $
\Ensure $\text{Cost Function in Equation \ref{eq:moment_optimization}} \leq \textbf{Tolerance}$
\While{$\lVert f(m(T)) \rVert \geq \textbf{Tolerance}$}
\State $m(0), m(\Delta t), \hdots, m(T) \gets F(t, m(t), U)$
\State $i \gets 1$
\While{$i \leq N_{steps}$}
\State $A_{i} \gets \frac{\partial}{\partial m}F((i)\Delta t, m((i)\Delta t), U_{i+1})$
\State $B_{i} \gets \frac{\partial}{\partial U}F((i)\Delta t, m((i)\Delta t), U_{i+1})$
\State $i \gets i+1$
\EndWhile
\State $H \gets [A_{N_{steps}}\hdots A_{2}B_{1} |A_{N_{steps-1}}\hdots A_{3}B_{2} | \hdots$ \par 
\hskip\algorithmicindent $| A_{N_{steps}}B_{N_{steps}-1} |B_{N_{steps}}]$
\State $\Delta U \gets \min_{\Delta U} \Delta U^{T} (H^{T}H+\lambda I)\Delta U$   \par 
\hskip\algorithmicindent $+ f(m(T))^{T}H\Delta U$
\State $U \gets U + \Delta U$
\EndWhile
\end{algorithmic}
\end{algorithm}
\vspace{-1ex}

The computational procedure used to synthesize ensemble controls is described in Algorithm \ref{alg:Iterative}, in which a cost function is iteratively reduced until a desired error tolerance is reached. The methodology can be described as consisting of three stages: (1) propagation of moment states using the current control profile; (2) linearization of time-localized dynamics; (3) and an update of the control profile. Both the propagation of moment states and respective linearizations are performed using a Taylor series expansion of higher order over the symbolic expression that represents the dynamics in Equation \eqref{eq:BEC_moment}. This approach is adopted for improved precision in the linearized approximation defined by the set of matrices $A_{k}$ and $B_{k}$, which are coupled in a matrix $H$, used to estimate the variation of the moment state at the final time step. The optimization problem can be represented using a quadratic program with a quadratic objective function and linear constraints, which we solve using a quadratic programming solver. This function has a penalty term scaled using a parameter $\lambda$, which suppresses variations in the defined control variation $\Delta U$, which would invalidate previously obtained linearized transformations. The related code is developed and executed on MATLAB™, version R2023a.

In our computational studies, the constraint bound values defined in Equation \eqref{eq:moment_optimization} are set to $\Omega_{\max} = 100$, $\Delta\Omega_{\min} = -500 s^{-1}$, and $\Delta\Omega_{\max} = 500 s^{-1}$. We examine control designs for various target states in momentum space with real controls where  $\Omega_{\min} = -100$, as well as with strictly non-negative real-valued control functions where $\Omega_{\min} = 0$.  We use a nominal value of $\Delta t = 0.001$ over a total time horizon of $T\in[0, 3]$,  resulting in 3000 variables and 12000 constraints for the quadratic program at each iteration.  The same momentum level truncation at $N^+=9$ and moment order $N=20$ are used in all computations.  We take an initial guess of $U_0(t)=2.5+sin(t)$. Our results are described in the following section. Table \ref{tbl:Simul_Times} contains a list of control synthesis problems specified by different values of parameters $n$ and $\Delta\Omega_{min}$, and the resulting time required to compute the ensemble controls.  We observe that the computation time is affected by the increase in complexity required to achieve a higher momentum level $n$. All simulations were performed on a 2023 MacBook Pro with Apple(R) M2 Pro Processor at 3.4 GHz and 16GB of memory.

\begin{table}[h!]
    \centering
    \caption{Computation times for moment pulse design.}
    \begin{tabular}{||c | c || c||}
     \hline
     $n$ in $C_f^+$ & $\Delta\Omega_{min}$ & Computation time (s) \\ [0.5ex] 
     \hline\hline
     1 & 0 & 866.4169\\ 
     \hline
     2 & 0 & 827.9643\\
     \hline
     3 & 0 & 1079.8969\\
     \hline
     4 & 0 & 1039.5165\\
     \hline\hline
     1 & -100 & 700.5696\\ 
     \hline
     2 & -100 & 845.3960\\
     \hline
     3 & -100 & 1001.5283\\
     \hline
     4 & -100 & 1088.8991\\
     \hline
    \end{tabular}
    \label{tbl:Simul_Times}
    \vspace{-2ex}
\end{table}


\begin{figure*}[t]
\begin{center}
\includegraphics[width=1.04\textwidth]{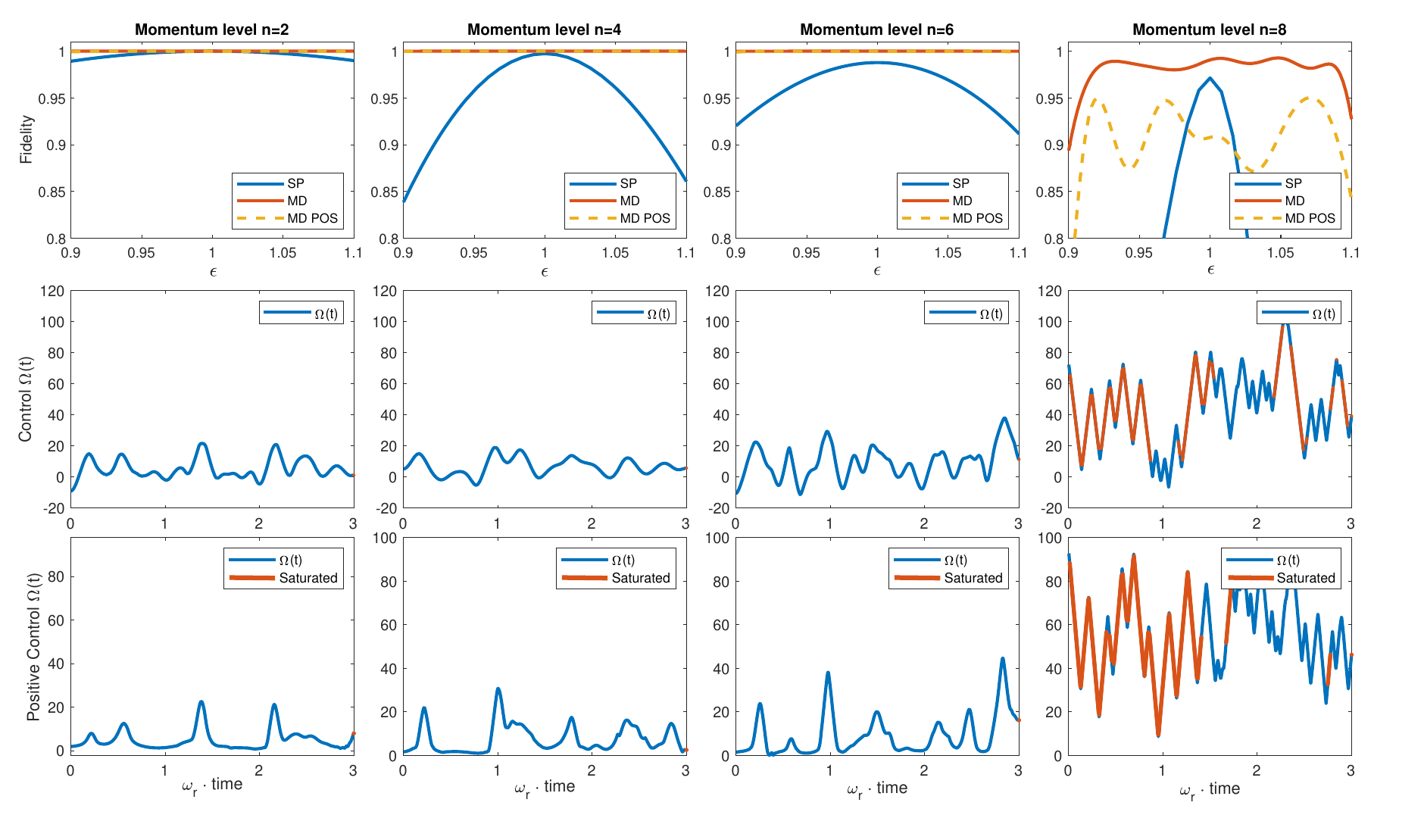}    
\vspace{-6ex}
\caption{Performance of moment ensemble controls compared to optimal square pulse sequences. Designs are for $\epsilon\in[1-\delta,1+\delta]$ with $\delta = 0.1$ and a time horizon of $T = 3$. Desired momentum levels are (from left to right) $n=2$, $n=4$, $n=6$, and $n=8$, corresponding to $C_{2}^{+}$, $C_{4}^{+}$, $C_{6}^{+}$, and $C_{8}^{+}$. The plots show (from top to bottom): terminal state fidelity $1-I_e$ with $I_e$ as in equation \eqref{eq:performance_index} for the square pulse (SP), the moment dynamics (MD) method, and MD method with positivity requirement (MD POS) as functions of the ensemble parameter $\epsilon$; the optimal control function without amplitude bounds; and optimal controls with strict non-negativity constraints. The times when the skew-rate limit $|\dot{\Omega}(t)|\leq 500$ is binding (Saturated) is indicated.}
\label{fig:Results_d01}
\end{center}
\vspace{-2ex}
\end{figure*}

\subsection{Moment Ensemble and Square Pulse Design} \label{ssec:MESPD}

To demonstrate the performance improvement of the Legendre moment dynamics (MD) method with respect to the square pulse (SP) sequence obtained by solving problem \eqref{eq:Square_Design_one}, we compare controls computed for an ensemble defined by $\delta = 0.1$. The square pulses are designed using a maximum momentum level $N^{+} = 24$, and the parameters that define the optimal pulses for target momentum states $n=1,2,3,4$ are given in Table \ref{tbl:Square_Pulse}. For control synthesis using the proposed MD ensemble method, we use a moment order $N = 20$ and truncate the momentum level at $N^{+} = 9$.  The resulting controls are validated by applying them to simulations of the entire ensemble with $\varepsilon\in[1-\delta,1+\delta]$ with the momentum levels truncated above $N^{+} = 24$. 

\begin{table}[h!]
    \centering
    \caption{Optimal square pulse design parameters.}
    \begin{tabular}{||c | c c c c c||} 
     \hline
     $n$ & $\Omega_{1}/\omega_{r}$ & $\Omega_{2}/\omega_{r}$ & $\omega_{r}\tau_{1}$ & $\omega_{r}\tau_{2}$ & $\omega_{r}\tau_{3}$ \\ [0.5ex] 
     \hline\hline
     1 & 3.9865 & 2.2849 & 0.4744 & 0.9427 & 0.4181\\ 
     \hline
     2 & 13.0036 & 9.5440 & 1.000 & 0.7190 & 1.000\\
     \hline
     3 & 32.4012 & 34.7591 & 0.1905 & 0.5523 & 0.1913 \\
     \hline
     4 & 41.4215 & 41.4263 & 1.2653 & 0.8002 & 1.9952 \\
     \hline
    \end{tabular}
    \label{tbl:Square_Pulse}
    \vspace{-2ex}
\end{table}

The real-valued and positive ensemble controls and terminal state error functions for the ensemble controls and square pulse sequences are compared in Figure \ref{fig:Results_d01} for target momentum levels $n = 1$, $2$, $3$, and $4$. Control fidelity shown there is defined as $1-I_e$, where $I_e$ is the performance index in equation \eqref{eq:performance_index}.  The square pulse sequence performs as expected, with error that is negligible at the nominal value of $\varepsilon = 1$, but which quickly increases as $\varepsilon$ diverges from unity.   In contrast, the terminal state error for the simulation using the controls obtained by solving problem \eqref{eq:moment_optimization} are applied to the ensemble remains quite low for all values of the ensemble parameter $\varepsilon\in[0.9,1.1]$ for target momentum levels $n = 1$, $2$, $3$, and still shows a significant performance advantage in the case of a target state $C_f^+=C_{f,2n}^+$ with $n=4$.  

The proposed method shows clear improvement with respect to all tested square pulse sequences in overall performance for momentum transfer of the ensemble of systems for positively-constrained controls as well as real-valued ones.  Observe that the amplitude and variation of the controls obtained for $n=4$ are greater than for the lower valued momentum state targets, which can be explained by the need to manipulate higher frequency dynamics.  Indeed, the constraints on amplitude and variation for control function are binding at many times during the optimization horizon, which necessarily limits the degree to which the algorithm can meet the optimization objective.  There are inherent trade-offs between the constraint bound values, the optimization horizon, and the complexity of the control task, which we have illustrated here.

\begin{figure*}
\centering
\vspace{-1ex}
\includegraphics[width=\textwidth]{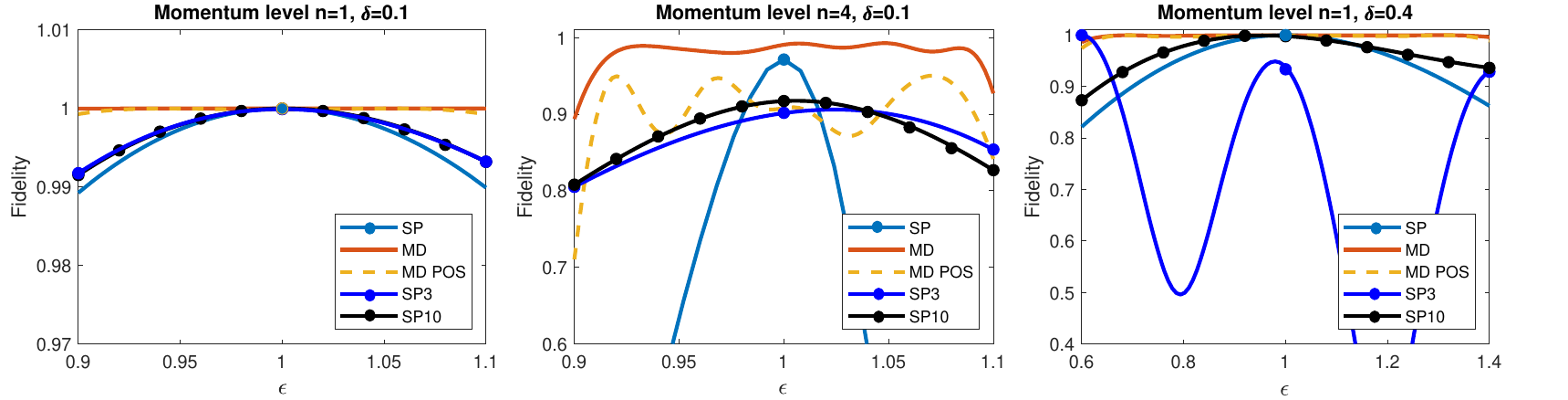}
\vspace{-3.5ex}
\caption{Control performance given as terminal fidelity $1-I_e$ with $I_e$ as in equation \eqref{eq:performance_index}, for the sample-design square pulse method for one (SP), three (SP3), and ten (SP10) samples in ensemble space, as well as the moment dynamics (MD) and positive MD method. The design parameters used are (from left to right): $\delta = 0.1$ and $n = 1$; $\delta = 0.1$ and $n = 4$; and $\delta = 0.4$ and $n = 1$.  Observe that the MD method clearly has the best performance.}
\label{fig:Results_All}
\vspace{-3ex}
\end{figure*}

\begin{table}[h!]
    \centering
    \caption{Performance results for moment ensemble controls for $\delta = 0.1$.}
    \begin{tabular}{||c | c c c||} 
     \hline
     $n$ & $I_{e, sp}$ & $I_{e, +}$ & $I_{e, \mathbb{R}}$\\ [0.5ex] 
     \hline\hline
     1 & 0.0141 & 0.0021 & 0.0004 \\ 
     \hline
     2 & 0.0579 & 0.0012 & 0.0004 \\
     \hline
     3 & 0.0515 & 0.0030 & 0.0006\\
     \hline
     4 & 0.1768 & 0.0608& 0.0260\\
     \hline
    \end{tabular}
    \vspace{-2ex}
    \label{tbl:Results_d01}
\end{table}
The performance achieved by the proposed method is quantified in Table \ref{tbl:Results_d01}, in which the index $I_{e}$ defined in Equation \eqref{eq:performance_index} is given for the square pulse controls and the positive and real valued ensemble controls as shown in Figure \ref{fig:Results_d01}, for which the index values are denoted by $I_{e, sp}$, $I_{e, +}$ and $I_{e, \mathbb{R}}$, respectively.  As seen in the top row of Figure \ref{fig:Results_d01}, the results in Table \ref{tbl:Results_d01} show that the positive valued ensemble controls result in a significant decrease in terminal error with respect to the square pulse sequence, and the real-valued ensemble controls are even better.  We note that in practice, the control represents the amplitude envelope of an optical pulse, so $\Omega(t)$ must remain positive in the present setting.  Reformulation of the RNE and adjustment of the experiment may enable real-valued control inputs, however.  

\subsection{Moment Ensemble and Ensemble Square Pulse Design}

We compare the ensemble controls to square pulse controls obtained by solving the optimization problem \eqref{eq:Square_Design_multi} for three control design cases: target $n = 1$ with $\delta = 0.1$; target $n = 4$ with $\delta=0.1$; and target $n = 1$ with $\delta = 0.4$.  The parameter values obtained by solving the square pulse design problem \eqref{eq:Square_Design_multi} using multiple $\varepsilon_i$ samples for $i=1,\ldots,m$ are given in Table \ref{tbl:Ens_Square_Pulse}. The results for these scenarios are shown in Figure \ref{fig:Results_All}, where circles indicate the error for the  values of $\epsilon$ sampled when solving problem \eqref{eq:Square_Design_multi} and lines show the error when applying the control to the continuum of values on the design interval $\epsilon\in[1-\delta,1+\delta]$.  The values of the performance index $I_e$ for the compared controls are shown in Table \ref{tbl:Results_Ens}.  The table includes values $I_{e,sp}$, $I_{e, +}$, $I_{e, \mathbb{R}}$, $I_{e, 3}$, and $I_{e, 10}$ for the nominal parameter square pulse, the positive ensemble control, the real-valued ensemble control, the square pulse with $m=3$ samples, and the square pulse with $m=10$ samples, respectively, for all three control design cases. 

\begin{table}[t!]
    \centering
    \caption{Ensemble Square pulse design parameters.}
    \begin{tabular}{||c c c | c c c c c||} 
     \hline
     $n$ & $\delta$ & $m$ & $\Omega_{1}/\omega_{r}$ & $\Omega_{2}/\omega_{r}$ & $\omega_{r}\tau_{1}$ & $\omega_{r}\tau_{2}$ & $\omega_{r}\tau_{3}$ \\ [0.5ex] 
     \hline\hline
     1 & 0.1 & 3 & 4.6404 & 4.2278 & 0.1943 & 0.9618 & 0.5286\\ 
     \hline
     1 & 0.4 & 3 & 10.6354 & 6.6525 & 1.2038 & 3.5000 & 0.5420\\ 
     \hline
     4 & 0.1 & 3 & 64.3286 & 27.5580 & 0.0924 & 3.4724 & 0.6647 \\
     \hline
     1 & 0.1 & 10 & 4.0446 & 4.2140 & 0.3769 & 0.6386 & 0.6767 \\
     \hline
     1 & 0.4 & 10 & 60.0000 & 4.6034 & 0.0159 & 2.6005 & 0.5322 \\
     \hline
     4 & 0.1 & 10 & 60.8008 & 30.2755 & 0.0967 & 3.4694 & 0.6798 \\
     \hline
    \end{tabular}
    \vspace{-5ex}
    \label{tbl:Ens_Square_Pulse}
\end{table}




\begin{table}[h!]
    \centering
    \caption{Performance including the ensemble square pulse approach.}
    \begin{tabular}{||c c | c c c c c||} 
     \hline
     $n$ & $\delta$ & $I_{e, sp}$ & $I_{e, +}$ & $I_{e, \mathbb{R}}$ & $I_{e, 3}$ & $I_{e, 10}$\\ [0.5ex] 
     \hline\hline
     1 & 0.1 & 0.0141 & 0.0021 & 0.0004 & 0.0127 & 0.0125\\ 
     \hline
     4 & 0.1 & 0.1768 & 0.0608 & 0.0260 & 0.0984 & 0.0965 \\
     \hline
     1 & 0.4 & 0.2278 & 0.0030 & 0.0006 & 0.5884 & 0.1818\\
     \hline
    \end{tabular}
    \label{tbl:Results_Ens}
    \vspace{-6ex}
\end{table}

We see in Figure \eqref{fig:Results_All} (right) that performance can be poor at values of $\epsilon$ within the design interval but far from the sample values (circles).  The results show the limitations of the square pulse design, even when multiple samples are used over the ensemble parameter domain. For the first case shown at left in Figure \ref{fig:Results_All}, the ensemble square pulse approach (Problem \eqref{eq:Square_Design_multi}) provides little improvement in relation to the nominal square pulse (Problem \eqref{eq:Square_Design_one}). The few degrees of freedom of the method and the small domain $K$ limits the advantage of sampling in the ensemble space. Although the advantage of sample-based design is noticeable in the second and third cases, the moment dynamics method shows major improvements over any square pulse designs.

We conclude that directly sampling the ensemble parameter space $K$ for the OCP does not effectively compensate for inhomogeneities in the experimental setting. In contrast, the proposed ensemble control method promises to achieve improved and homogeneous performance for an entire BEC ensemble in a matter wave interferometry experiment.




\section{CONCLUSION} \label{sec:conc}


We present a method for designing constrained optical pulses for optimal matter-wave splitting of a Bose-Einstein condensate in the presence of experimental inhomogeneities, which induce a parameter uncertainty within the momentum evolution of the collection of atoms. An ensemble of systems is parameterized by the factor that modulates the optical pulse, and Legendre moments are used to represent the optimal control problem for an infinite dimensional system as an optimization problem in moment space. The moment dynamics approach effectively compensates for inhomogeneities on a continuum by using a spectral approximation of the ensemble in Hilbert space.  The spectral accuracy of polynomial approximation on an orthogonal basis ensures exponential convergence to the continuum state with increasing approximation order. This representation enables design of control sequences for ensemble systems in a space of greatly reduced dimension without significant loss of precision \cite{Li:22}.  Crucially, spectral methods such as the one developed here provide guarantees for the control design for all values in the entire uncertain parameter space, in contrast to sampling approaches usually used in stochastic optimization that, e.g., do not provide guarantees for values in intervals between samples. Convergence analysis of iterative methods for bilinear systems remains an active area of research \cite{wang2017fixed}.

The fidelity achieved by the proposed method has advantages over recently developed state-of-the-art optical pulse sequences.  In particular, the method achieves precise momentum transfer of the BEC to high diffraction orders with inherent robustness to inhomogeneity in the effect of the optical pulse on atoms in the sample, and yields continuously-varying pulses that can be tuned to maximize the effectiveness of equipment in experimental settings.
We expect future research to extend our results to more detailed models of BEC-splitting, as well as control protocols for additional procedures for state preparation and measurement.  Moment dynamics can be extended to incorporate more parameters to improve robustness to additional sources of uncertainty and inhomogeneity.  Our approach could be improved and analyzed to characterize the trade-offs between truncations in the moment and quantum Hilbert spaces, constraint bound values, pulse duration, and control performance.  Based on the successful validation of our approach in simulation, we expect that optimal pulses can be tested in experimental settings to confirm anticipated improvements in performance of matter-wave diffraction metrology. 

\vspace{2ex}

\addtolength{\textheight}{-12cm}   




\bibliographystyle{unsrt}  

\typeout{}
\bibliography{refs}

\end{document}